
\input amstex.tex
\loadeusm
\loadeusb
\documentstyle{amsppt}

\topmatter

\title
The Cofinality Spectrum Of The Infinite Symmetric Group
\endtitle

\author
Saharon Shelah\\
and\\
Simon Thomas
\endauthor
\address
Department of Mathematics\\
Rutgers University\\
New Brunswick, New Jersey 08903\\
USA
\endaddress
\thanks Research partially supported by the BSF. Publication 524 of
the first author.
\medskip
Research partially supported by NSF Grants.
\endthanks

\endtopmatter
\newpage
\subhead 1.  Introduction\endsubhead
\medskip
Suppose that $G$ is a group that is not finitely generated.  Then $G$
can be written as the union of a chain of proper subgroups.  The {\it
cofinality spectrum} of $G$, written $CF(S)$, is the set of regular
cardinals $\lambda$ such that $G$ can be expressed as the union of a
chain of $\lambda$ proper subgroups.  The {\it cofinality} of $G$,
written $c(G)$, is the least element of $CF(G)$.
\medskip
Throughout this paper, $S$ will denote the group Sym$(\omega)$ of all
permutations of the set of natural numbers.  In [MN], Macpherson and
Neumann proved that $c(S) > \aleph_o$. In [ST1] and [ST2], the
possibilities for the value of $c(S)$ were studied.  In particular, it
was shown that it
is consistent that $c(S)$ and $2^{\aleph_o}$ can be any two prescribed
regular uncountable cardinals, subject only to the obvious requirement
that $c(S) \leq 2^{\aleph_o}$.  In this paper, we shall
begin the study of the possibilities for the set $CF(S)$.
\medskip
There is one obvious constraint on the set $CF(S)$,
arising from the fact that $S$ can be expressed as the union of a
chain of $2^{\aleph_o}$ proper subgroups; namely, that
$cf(2^{\aleph_o}) \in CF(S)$. 
  Initially it is difficult to think of any
other constraints on $CF(S)$.  And we shall show that it is consistent
that $CF(S)$ is quite a bizarre set of cardinals.  For example, the
following result is a special case of our main theorem.
\proclaim{Theorem 1.1}
Let $T$ be any subset of $\omega \smallsetminus \{0\}$.  Then it is consistent
that\linebreak $\aleph_n \in CF(S)$ if and only if $n\in T$.\endproclaim
\medskip
After seeing this result, the reader might suspect that it is
consistent that $CF(S)$ is an arbitrarily prescribed set of regular
uncountable cardinals, subject only to the above mentioned constraint.
However, this is not the case.
\proclaim{Theorem 1.2}
If $\aleph_n \in CF(S)$ for all $n\in
\omega \smallsetminus \{0\}$, then $\aleph_{\omega +1} \in CF(S)$.
\endproclaim\medskip
(Of course, this result is only interesting when $2^{\aleph_0} >
\aleph_{\omega +1}$.)  In Section 2, we shall use $pcf$ theory to
prove Theorem 1.2, together with some further results which restrict
the possibilities for $CF(S)$.  In Section 3, we shall prove the
following result.
\proclaim{Theorem 1.3}
Suppose that $V \vDash GCH$.  Let $C$ be a set of regular uncountable
cardinals which satisfies the following conditions.
\medskip
(1.4)
\roster
\item"{(a)}"  $C$ contains a maximum element.
\item"{(b)}"  If $\mu$ is an inaccessible cardinal such that $\mu=\sup
(C\cap\mu)$, then $\mu \in C$.
\item"{(c)}"  If $\mu$ is a singular cardinal such that $\mu = \sup(C\cap\mu)$, then $\mu^+ \in C$.
\endroster
Then there exists a $c.c.c$ notion of forcing $\Bbb P$ such that
$V^{\Bbb P} \vDash CF(S)=C$. 
\endproclaim
\medskip
This is not the best possible result. In particular, clause (1.4)(c) can be improved
 so that we gain a little more control over what
occurs at successors of singular cardinals.  This matter will be
discussed more fully at the end of Section 2. Also clause (1.4)(a) is
not a necessary condition. For example, let $V \vDash GCH$ and let
$C=\{\aleph_{\alpha +1} \alpha < \omega_1\}$. At the end of Section 3,
we shall show that if $\kappa$ is any singular cardinal such that
$cf(\kappa) \in C$, then there exists a $c.c.c$ notion of forcing
$\Bbb P$ such that
$V^{\Bbb P} \vDash CF(S)=C \text{ and } 2^{\aleph_o} = \kappa$. In particular,
$2^{\aleph_o}$ cannot be bounded in terms of the set $CF(S)$.
\medskip
In this paper, we have made no attempt to control what occurs at inaccessible
cardinals $\mu$ such that $\mu = \sup(C\cap\mu)$. We intend to deal with this
matter in a second paper, which is in preparation. In this second paper, we
also hope to give a complete characterisation of those sets $C$ for which there
exists a $c.c.c$ notion of forcing $\Bbb P$ such that 
$V^{\Bbb P} \vDash CF(S)=C$.
\medskip
Our notation mainly follows that of Kunen [K].  Thus if $\Bbb P$ is a
notion of forcing and $p,q\in \Bbb P$, then $q\leq p$ means that $q$
is a strengthening of $p$.  If $V$ is the ground model, then we often
denote the generic extension by $V^\Bbb P$ if we do not wish to
specify a particular generic filter $G\subseteq \Bbb P$.  If we want
to emphasize that the term $t$ is to be interpreted in the model $M$
of $ZFC$, then we write $t^M$; for example, Sym$(\omega)^M$.  If
$A\subseteq \omega$, then $S_{(A)}$ denotes the pointwise stabilizer
of $A$.  Fin$(\omega)$ denotes the subgroup of elements $\pi\in S$
such that the set $\{n<\omega \pi(n)\neq n\}$ is finite.  If $\phi,
\psi \in S$, then we define $\phi =^\ast \psi$ if and only if
$\phi\psi^{-1} \in$ Fin$(\omega)$.

\subhead 2. Some applications of $pcf$ theory\endsubhead
\medskip
Let $\langle \lambda_i i\in I\rangle$ be an indexed set of regular
cardinals.  Then $\underset{i\in I}\to{\Pi}\lambda_i$ denotes the set
of all functions  $f$ such that dom $f=I$ and $f(i)\in \lambda_i$ for
all $i\in I$.  If $\Cal F$ is a filter on $I$ and $\Cal I$ is the dual
ideal, then we write either $\underset{i\in I}\to{\Pi} \lambda_i/_\Cal
F$ or $\underset{i\in I}\to{\Pi}\lambda_i/_\Cal I$ for the
corresponding reduced product.  We shall usually prefer to work with
functions $f\in \underset{i\in I}\to{\Pi}\lambda_i$ rather than with
the corresponding equivalence classes in $\underset{i\in I}\to {\Pi}
\lambda _i/_\Cal I$.  For $f,g,\in \underset{i\in I}\to{\Pi}\lambda_i$, we
define
$$
f\leq_{\Cal I} g \text{ iff }\{i\in If(i) > g(i) \} \in \Cal I
$$
$$
f <_{\Cal I} g \text{ iff } \{ i\in If(i) \geq g(i)\}\in \Cal I.
$$
We shall sometimes write $f\leq_{\Cal F} g, f<_{\Cal F} g$ instead of
$f\leq_{\Cal I} g, f<_{\Cal I} g$ respectively.  If  $\Cal I=\{\phi\}$,
then we shall write $f\leq g, f<g$.  Suppose that there exists a
regular cardinal $\lambda$ and a sequence $\langle
f_\alpha|\alpha<\lambda\rangle$ of elements of $\underset{i\in
I}\to{\Pi}\lambda_i$ such that 
\roster
\item"{(a)}"  if $\alpha < \beta < \lambda$, then $f_\alpha <_{\Cal
I}f_\beta$; and 
\item"{(b)}"  for all $h\in \underset{i\in I}\to {\Pi}\lambda_i$, there exists
$\alpha < \lambda$ such that $h<_{\Cal I} f_\alpha$.
\endroster
Then we say that $\lambda$ is the {\it true cofinality} of
$\underset{i\in I}\to{\Pi}\lambda_i/_\Cal I$, and write $tcf\left(\underset{i\in
I}\to{\Pi} \lambda_i/_{\Cal I}\right)= \lambda$.
Furthermore, we say that $\langle f_\alpha|\alpha < \lambda\rangle$
{\it witnesses} that $tcf\left(\underset{i\in
I}\to{\Pi}\lambda_i/_{\Cal I}\right)=\lambda$. For example, if $\Cal D$
is an ultrafilter on $I$, then $\underset{i\in I}\to{\Pi}
\lambda_i/_{\Cal D}$ is a linearly ordered set and hence has a true
cofinality.  A cardinal $\lambda$ is a {\it possible cofinality} of
$\underset{i\in I}\to{\Pi}\lambda_i$ if there exists an ultrafilter
$\Cal D$ on $I$ such that $tcf\left(\underset{i\in I}\to{\Pi}\lambda_i /_{\Cal
D}\right) = \lambda$.  The set of all possible cofinalities of
$\underset{i\in I}\to{\Pi}\lambda_i$ is $pcf\left(\underset{i\in I}\to{\Pi}
\lambda_i\right)$.
\medskip
In recent years, Shelah has developed a deep and beautiful theory of
the structure of $pcf\left(\underset{i\in I}\to{\Pi}\lambda_i\right)$
when $|I|<\min\{\lambda_i|i\in I\}$.  A thorough development of $pcf$
theory and an account of many of its applications can be found in [Sh-g].
[BM] is a self-contained survey of the basic elements of $pcf$ theory.
In this section of the paper, we shall see that $pcf$ theory imposes a
number of constraints on the possible structure of $CF(S)$.  (Whenever
it is possible, we shall give references to both [Sh-g] and [BM] for
the results in $pcf$ theory that we use.)
\proclaim{Theorem 2.1} 
Suppose that $\langle \lambda_n|n<\omega\rangle$ is a strictly
increasing sequence of cardinals such that $\lambda_n \in CF(S)$ for
all $n<\omega$.  Let $\Cal D$ be a nonprincipal ultrafilter on
$\omega$, and let $tcf\left(\underset{n<\omega}\to{\Pi}
\lambda_n/_{\Cal D}\right)=\lambda$.  Then $\lambda \in CF(S)$.
\endproclaim
\demo{Proof}
For each $n<\omega$, express $S=\underset{i<\lambda_n}\to{\bigcup} G^n_i$
as the union of a chain of $\lambda_n$ proper subgroups.  Let
$\langle f_\alpha|\alpha <\lambda\rangle$ be a sequence in
$\underset{n<\omega}\to{\Pi}\lambda_n$ which witnesses that
$tcf\left(\underset{n<\omega}\to{\Pi}\lambda_n/_{\Cal
D}\right)=\lambda$.  For each $\alpha < \lambda$, let $H_\alpha$ be
the {\it set} of all $g\in S$ such that $\{n<\omega|g \in
G^n_{f_\alpha(n)}\} \in \Cal D$.  Then it is easily checked that
$H_\alpha$ is a subgroup of $S$, and that $H_\alpha \subseteq H_\beta$
for all $\alpha < \beta < \lambda$.  Suppose that $g\in S$ is an
arbitrary element.  Define $f\in \underset{n <\omega}\to{\Pi}
\lambda_n$ by $f(n)= \min\{i|g\in G^n_i\}$.  Then there exists $\alpha
< \lambda$ such that $f<_{\Cal D} f_\alpha$.  Hence $g\in H_\alpha$.
Thus $S=\underset{\alpha <\lambda}\to{\bigcup} H_\alpha$.
\medskip 
So it suffices to prove that $H_\alpha$ is a proper subgroup of $S$
for each $\alpha < \lambda$.  Fix some $\alpha < \lambda$.  Lemma 2.4
[MN] implies that for each $n<\omega$, $i<\lambda_n$ and $X\in
[\omega]^\omega$, the setwise stabilizer of $X$ in $G^n_i$ does not
induce Sym $(X)$ on $X$.  Express $\omega =
\underset{n<\omega}\to{\bigcup} X_n$ as the disjoint union of countably
many infinite subsets $X_n$.  For each $n<\omega$, choose $\pi_n \in 
$ Sym $(X_n)$ such that $g\upharpoonright X_n\neq \pi_n$ for all
$g\in G^n_{f_\alpha(n)}$.  Then $\pi=\underset{n<\omega}\to{\bigcup} \pi_n
\in S\smallsetminus H_\alpha$.
\enddemo
$\quad \hfill\qed$

\demo{Proof of Theorem 1.2}
\medskip
By [Sh-g, II 1.5] (or see [BM, 2.1]), there exists an ultrafilter
$\Cal D$ on $\omega$ such that $tcf\left(\underset{n<\omega}\to{\Pi}
\aleph_n/_{\Cal D}\right) = \aleph_{\omega +1}$.
\enddemo
$\quad\hfill\qed$
\medskip
If we assume $MA_\kappa$, then we can obtain the analogous result for
cardinals $\kappa$ such that $\aleph_o < \kappa < 2^{\aleph_o}$.  (In Section 3, we
shall prove that the following result cannot be proved in $ZFC$.)
\proclaim{Theorem 2.2 ($MA_\kappa$)}
Suppose that $\langle \lambda_\alpha|\alpha < \kappa\rangle$ is a strictly
increasing sequence of cardinals such that $\lambda_\alpha \in CF(S)$
for all $\alpha < \kappa$.  Let $\Cal D$ be a nonprincipal ultrafilter on
$\kappa$, and let $tcf\left(\underset{\alpha <
\kappa}\to{\Pi}\lambda_\alpha/_{\Cal D}\right)=\lambda$. Then
$\lambda\in CF(S)$.
\endproclaim
\demo{Proof} For each $\alpha < \kappa$, express
$S=\underset{i<\lambda_\alpha}\to{\bigcup} G^\alpha_i$ as the union of a
chain of $\lambda_\alpha$ proper subgroups.  Let $\langle f_\beta
|\beta <\lambda\rangle$ be a sequence in $\underset{\alpha <
\kappa}\to{\Pi} \lambda_\alpha$ which witnesses that
$tcf\left(\underset{\alpha < \kappa}\to{\Pi}\lambda_\alpha/_{\Cal
D}\right)=\lambda$.  For each $\beta < \lambda$, let $H_\beta$ be the
set of all $g\in S$ such that $\{\alpha <\kappa |g\in
G^\alpha_{f_\beta(\alpha)}\}\in \Cal D$.  Arguing as in the proof of
Theorem 2.1, it is easily checked that $\langle H_\beta |\beta <
\lambda\rangle$ is a chain of subgroups such that $S= \underset{\beta
< \lambda}\to{\bigcup} H_\beta$.
\medskip 
Thus it suffices to prove that $H_\beta$ is a proper subgroup of $S$
for each $\beta < \lambda$.  Fix some $\beta <\lambda$.  Suppose that
we can find an element $g\in S\smallsetminus\underset{\alpha <
\kappa}\to{\bigcup} G^\alpha_{f_\beta(\alpha)}$.  
\medskip
Then clearly $g\notin H_\beta$.  But the existence of such an element
$g$ is an immediate consequence of the following theorem.
\enddemo
$\quad\hfill\qed$
\proclaim{Theorem 2.3 ($MA_\kappa$)} Suppose that for each $\alpha < \kappa,
S=\underset{i<\theta_\alpha}\to{\bigcup} H^\alpha_i$ is the union of the
chain of proper subgroups $H^\alpha_i$.  Then for each
$f\in\underset{\alpha < \kappa}\to{\Pi} \theta_\alpha,
S\neq\underset{\alpha<\kappa}\to{\bigcup} H^\alpha_{f(\alpha)}$.
\endproclaim
\remark{Remark 2.4}  In [ST 1], it was shown that $MA_\kappa$ implies
that $c(S)>\kappa$.  This result is an easy consequence of Theorem 2.3.
\endremark
\remark{Remark 2.5}  In [MN], Macpherson and Neumann proved that if
$\{H_n|n<\omega\}$ is an {\it arbitrary} set of proper subgroups of
$S$, then $S\neq\underset{n<\omega}\to{\bigcup} H_n$.  It is an open
question whether $MA_\kappa$ implies the analogous statement for
cardinals $\kappa$ such that\linebreak $\aleph_o <\kappa < 2^{\aleph_o}$.
Regard $S$ as a Polish space in the usual way.  Then the proof of
Theorem 2.3 shows that the following result holds.
\endremark
\proclaim{Theorem 2.6 $(MA_\kappa)$}  Suppose that for each $\alpha <
\kappa, H_\alpha$ is a nonmeagre proper subgroup of $S$.  Then $S\neq
\underset{\alpha < \kappa}\to{\bigcup} H_\alpha$.  
\endproclaim
$\quad\hfill\qed$
\medskip
Unfortunately there exist maximal subgroups $H$ of $S$ such that $H$
is meagre.  For example, let $\omega=\Omega_1\cup \Omega_2$ be a
partition of $\omega$ into two infinite pieces.  Let 
$$
H=\{g\in S\big\vert |g[\Omega_1] \triangle \Omega_i |<\aleph_o \text{ for
some } i\in \{1,2\}\}.
$$     
(Here $\triangle$ denotes the symmetric difference.)  Then $H$ is a
maximal subgroup of $S$; and it is easily checked that $H$ is meagre.
\demo{Proof of Theorem 2.3 $(MA_\kappa)$}
We shall make use of the technique of generic sequences of elements of
$S$, as developed in [HHLSh]. (The slight differences in notation
between this paper and [HHLSh] arise from the fact that permutations
act on the left in this paper.)
\enddemo

\definition{Definition 2.7}  A finite sequence $\langle
g_1,\cdots,g_n\rangle \in S^n$ is {\it generic} if the following two
conditions hold.
\roster
\item  For all $A\in [\omega]^{<\omega}$, there exists $A\subseteq B \in
[\omega]^{<\omega}$ such that $g_i[B]=B$ for all $1\leq i\leq n$.
\item  Suppose that $A\in [\omega]^{<\omega}$ and that $g_i[A]=A$ for
all $1\leq i\leq n$.  Suppose further that $A\subseteq B\in
[\omega]^{<\omega}$ and that $h_i\in Sym(B)$ extends
$g_i\upharpoonright A$ for all $1\leq i\leq n$.  Then there exists
$\pi \in S_{(A)}$ such that $\pi g_i \pi^{-1}$ extends $h_i$ for all
$1\leq i\leq n$.
\endroster
\enddefinition
\proclaim{Claim 2.8}  If $\langle g_1,...,g_n\rangle,\langle
h_1,...,h_n\rangle \in S^n$ are generic, then there exists $f\in S$
such that $fg_if^{-1}=h_i$ for all $1\leq i\leq n$.
\endproclaim
\demo{Proof of Claim 2.8} This follows from Proposition 2.3 [HHLSh].
\enddemo
 $\quad\hfill\qed$
\medskip
>From now on, regard $S$ as a Polish space in the usual way.
\proclaim{Claim 2.9}  The set $\{\langle g_1,...,g_n\rangle \in
S^n\big\vert\langle g_1,...,g_n\rangle$ is generic$\}$ is comeagre in
$S^n$ in the product topology.
\endproclaim
\demo{Proof of Claim 2.9}  This follows from Theorem 2.9 [HHLSh].
\enddemo
$\quad\hfill\qed$
\proclaim{Claim 2.10}  If $\langle g_1,...,g_{n+1}\rangle \in S^{n+1}$
is generic, then for each $A\in [\omega]^{<\omega}, m\in
\omega\smallsetminus A$ and $1\leq \ell\leq n+1$, the following
condition holds.
\newline
$(2.11)_{A,m,\ell}$  Let $\Omega=\{i\big\vert 1\leq i\leq n+1,i\neq
\ell\}$.  If $g_i[A]=A$ for all $i\in \Omega$, then there exists $B\in
[\omega\smallsetminus A]^{<\omega}$ such that
\roster
\item"{(a)}"  $m\in B$;
\item"{(b)}"  $g_i[B]=B$ for all $i\in \Omega$;
\item"{(c)}"  $g_\ell[A\cup B]=A\cup B$;
\item"{(d)}"  for all $\pi\in$ Sym $(\Omega)$, there exists $\phi \in$
Sym $(B)$ such that $\phi(g_i\upharpoonright
B)\phi^{-1}=g_{\pi(i)}\upharpoonright B$ for all $i\in \Omega$.
\endroster
\endproclaim
\demo{Proof of Claim 2.10}  For each $A\in [\omega]^{<\omega}, m\in
\omega \smallsetminus A$ and $1\leq \ell \leq n + 1$, \newline
let
$C^{n+1}(A,m,\ell)$ consist of the sequences
$\langle g_1,...,g_{n+1}\rangle\in S^{n+1}$ which satisfy\newline
$(2.11)_{A,m,\ell}$.  Then it is easily checked that
$C^{n+1}(A,m,\ell)$ is open and dense in $S^{n+1}$.  Hence
$C^{n+1}=\underset{A,m,\ell}\to{\bigcap}C^{n+1}(A,m,\ell)$ is
comeagre in $S^{n+1}$.  Claim 2.9 implies that there exists a generic
sequence $\langle g_1,...,g_{n+1}\rangle \in C^{n+1}$.  So the result
follows easily from Claim 2.8.
\enddemo
$\quad\hfill\qed$
\definition{Definition 2.12}  If $\sigma$ is an infinite ordinal, then
the sequence $\langle g_i\big\vert i < \sigma\rangle$ of elements of
$S$ is {\it generic} if for every finite subsequence $ i_1< ... < i_n <
\sigma,\langle g_{i_1},...,g_{i_n}\rangle$ is generic.
\enddefinition
\medskip
We have now developed enough of the theory of generic sequences to
allow us to begin the proof of Theorem 2.3.  Consider the chains of
proper subgroups,\linebreak $S=\underset{i<\theta_\alpha}\to{\bigcup}
H^\alpha_i$ for $\alpha < \kappa$.  We can assume that Fin$(\omega)\leq
H^\alpha_o$ for all $\alpha < \kappa$.  Let\linebreak $f\in\underset{\alpha
<\kappa}\to{\Pi} \theta_\alpha$.  We must find an element $\pi\in
S\smallsetminus\underset{\alpha <\kappa}\to{\bigcup}
H^\alpha_{f(\alpha)}$.  We shall begin by inductively constructing a
generic sequence of elements of $S$.
$$
\langle g^o_o, g^1_o,...,g^o_\alpha, g^1_\alpha,...\rangle_{\alpha
<\kappa}
$$
such that for all $\alpha < \kappa$, there exist $f(\alpha)\leq
\gamma_\alpha <\theta_\alpha$ such that $g^o_\alpha \in
H^\alpha_{\gamma_\alpha}$ and $g^1_\alpha \notin
H^\alpha_{\gamma_\alpha}$.  Then we shall find an element $\pi\in S$
such that $\pi g^o_\alpha \pi^{-1} =^\ast g^1_\alpha$ for all $\alpha
<\kappa$.  This implies that $\pi \notin \underset{\alpha <
\kappa}\to{\bigcup} H^\alpha_{\gamma_\alpha} \supseteq
\underset{\alpha <\kappa}\to{\bigcup} H^\alpha_{f(\alpha)}$.
\medskip
Suppose that we have constructed $g^o_\beta, g^1_\beta$ for $\beta <
\alpha$.  For each finite subsequence $\bar g$ of $\langle g^o_\beta,
g^1_\beta \big\vert \beta < \alpha\rangle$, the set $\{h\in
S\big\vert\bar g\sphat\,  h\text{ is generic} \}$ is comeagre in $S$.
(See[HHLSh], page 216.)  Since $MA_\kappa$ implies that the union of
$\kappa$ meagre subsets of a Polish space is meagre, the set
$$
\{h \in S\big\vert \langle g^o_\beta, g^1_\beta |\beta <\alpha\rangle\sphat\, h \text{ is generic }\}
$$    
is also comeagre in $S$.  So we can choose a suitable $g^o_\alpha$ and
$f(\alpha)\leq \gamma_\alpha < \theta_\alpha$ with $g^o_\alpha \in
H^\alpha_{\gamma_\alpha}$.  The set
$$
C =\{h\in S\big\vert \langle g^o_\beta, g^1_\beta\beta <\alpha
\rangle \sphat\, g^o_{\alpha} {\sphat} \,h \text{ is generic} \}
$$
is also comeagre in $S$.  Since $H^\alpha_{\gamma_\alpha}$ is a proper
subgroup of $S$, we have that\linebreak $C\smallsetminus H^\alpha_{\gamma_\alpha} \neq \emptyset$.
(If not, then $H^\alpha_{\gamma_\alpha}$ is comeagre and hence so are
each of its cosets in $S$.  As any two comeagre subsets of $S$
intersect, this is impossible.)  Hence we can choose a suitable
$g^1_\alpha \in C\smallsetminus H^\alpha_{\gamma_\alpha}$.  Thus the
desired generic sequence can be constructed.
\proclaim{Lemma 2.13}  Let $\langle g^o_\alpha, g^1_\alpha \alpha
<\kappa\rangle$ be a generic sequence of elements of $S$.  Then there
exists a $\sigma$-centred notion of forcing $\Bbb P$ such that
$$
\underset{\Bbb P}\to{\Vdash}  \text{There exists } \pi\in \text{ Sym }(\omega)
\text{ such that } \pi g^o_\alpha \pi^{-1}=^\ast g^1_\alpha \text{ for
all } \alpha <\kappa.
$$
\endproclaim
\demo{Proof of Lemma 2.13}  Let $\Bbb P$ consist of the conditions
$p=\langle h, F\rangle$ such that 
\roster
\item"{(a)}" there exists $A\in [\omega]^{<\omega}$ such that $h\in$
Sym$(A)$;
\item"{(b)}"  $F\in [\kappa]^{<\omega}$;
\item"{(c)}"  for each $\alpha \in F$ and $\tau\in\{0,1\},\, 
g^\tau_\alpha[A]=A$.
\endroster
We define $\langle h_2,F_2\rangle\leq\langle h_1,F_1\rangle$ iff the
following two conditions hold.
\medskip
(1)  $h_1\subseteq h_2$ and $F_1 \subseteq F_2$.
\medskip
(2)  Let $B=$ dom $h_2\smallsetminus$ dom $h_1$ and let
$\phi=h_2\upharpoonright B$.  Then $\phi(g^o_\alpha\upharpoonright
B)\phi^{-1} = g^1_\alpha \upharpoonright B$ for each $\alpha \in F_1$.
\medskip
Clearly $\Bbb P$ is $\sigma-$centered.  Claim 2.10 implies that each
of the sets
$$
D_m =\{\langle h, F\rangle \big\vert m \in \text{ dom } h\} \quad
,m<\omega
$$
and 
$$
E_\alpha = \{\langle h,F\rangle \big\vert\alpha \in F\} \quad , \alpha <
\kappa,
$$
are dense in $\Bbb P$.  The result follows.
$\qquad\hfill\qed$
\medskip
This completes the proof of Theorem 2.3.
$\qquad\hfill\qed$
\medskip
The following theorem goes some way towards explaining why we have
assumed that $C$ satisfies condition (1.4)(c) in the statement of
Theorem 1.3.  (We will discuss this matter fully after we have proved
Theorem 2.15.)
\definition{Definition 2.14}  If $\delta$ is a limit ordinal, then
$J^{bd}_\delta$ is the ideal on $\delta$ defined by
$$  
J^{bd}_\delta = \{ B\big\vert \text{There exists } i < \delta
\text{ such that } B\subseteq i\}.
$$
\enddefinition
\proclaim{Theorem 2.15}  Let $\kappa$ be a regular cardinal, and
suppose that $\langle \lambda_\alpha |\alpha <\kappa\rangle$ is a
strictly increasing sequence of cardinals such that $\lambda_\alpha
\in CF(S)$ for all $\alpha < \kappa$.  Suppose further that $tcf\left(\underset{\alpha < \kappa}\to{\Pi}
\lambda_\alpha/_{J^{bd}_\kappa}\right) = \lambda$.
Then either $\kappa \in CF(S)$ or $\lambda \in CF(S)$.
\endproclaim
\demo{Proof}  Suppose that $\kappa \notin CF(S)$.  For each $\alpha <
\kappa$, express $S=\underset{i<\lambda_\alpha}\to{\bigcup}
G^\alpha_i$ as the union of a chain of $\lambda_\alpha$ proper
subgroups.  Let $\langle f_\beta|\beta <\lambda\rangle$ be a sequence
in $\underset{\alpha < \kappa}\to{\Pi} \lambda_\alpha$ which witnesses
that
$tcf\left(\underset{\alpha<\kappa}\to{\Pi}\lambda_\alpha/_{J^{bd}_\kappa}\right)=\lambda$.
For each $\beta <\lambda$, let $G^\ast_\beta$ be the set of all $g\in
S$ such that $\kappa\smallsetminus\{\alpha < \kappa|g\in
G^\alpha_{f_\beta(\alpha)}\} \in J^{bd}_\kappa$.  Arguing as before,
it is easily checked that $\langle G^\ast_\beta|\beta <
\lambda\rangle$ is a chain of subgroups such that $S=\underset{\beta <
\lambda}\to{\bigcup} G^\ast_\beta$.
\medskip
Thus it suffices to prove that $G^\ast_\beta$ is a proper subgroup of
$S$ for each $\beta <\lambda$.  So suppose that $G^\ast_\beta = S$ for
some $\beta < \lambda$.  For each $\alpha < \kappa$, define\newline
$H_\alpha=\bigcap\{G^\gamma_{f_\beta(\gamma)}|\alpha \leq \gamma <
\kappa\}$.  Then $\langle H_\alpha|\alpha < \kappa\rangle$ is a chain
of subgroups such that $S=\underset{\alpha < \kappa}\to{\bigcup}
H_\alpha$.  If $\alpha < \kappa$, then $H_\alpha \leq
G^\alpha_{f_\beta(\alpha)}$ and so $H_\alpha$ is a proper subgroup of
$S$.  But this contradicts the assumption that $\kappa\notin CF(S)$.
\enddemo
$\quad\hfill\qed$
\medskip
Suppose that $V\vDash GCH$, and that $\mu$ is a singular cardinal.
Let $\langle\theta_i|i<\eta\rangle$ be the strictly increasing
enumeration of all regular uncountable cardinals $\theta$ such that
$\theta < \mu$.  Let $\Cal F = \underset{i<\eta}\to{\Pi} \theta_i$.
Then $|\Cal F|=\mu^+$.  Now let $\Bbb P$ be any c.c.c. notion of
forcing.  From now on, we shall work in $V^\Bbb P$.  Since $\Bbb P$ is
c.c.c., for each $g\in \underset{i<\eta}\to{\Pi} \theta_i$, there
exists $f\in \Cal F$ such that $g\leq f$.  Suppose now that
$\langle\lambda_\alpha|\alpha < \delta\rangle$ is an increasing
subsequence of $\langle\theta_i|i<\eta\rangle$ such that
$|\delta|<\lambda_o$ and $\underset{\alpha<\delta}\to{\sup}\, \lambda_\alpha =\mu$.  Let
$$
\Cal F^\ast =\{f\in\underset{\alpha < \delta}\to{\Pi}\lambda_\alpha
\big\vert \text{ There exists } h\in \Cal F \text{ such that }
f\subseteq h\}.
$$  
\medskip
Then for all $g\in \underset{\alpha<\delta}\to{\Pi}\lambda_\alpha$,
there exists $f\in \Cal F^\ast$ such that $g\leq f$. This implies
that $\max (pcf(\underset{\alpha <
\delta}\to{\Pi}\lambda_\alpha))=\mu^+$.  By [Sh-g,I] (or see [BM,4.3]),
we obtain that
$tcf(\underset{\alpha<\delta}\to{\Pi}\lambda_\alpha/_{J^{bd}_\delta})
= \mu^+$.  In summary, we have shown that the following statement is
true in $V^{\Bbb P}$.
\medskip
{\bf The Strong Hypothesis (2.16).}
Let $\delta$ be a limit ordinal, and let $\langle\lambda_\alpha|\alpha
<\delta\rangle$ be a strictly increasing sequence of regular cardinals
such that $|\delta<\lambda_o$.  Then $tcf\left(\underset{\alpha <
\delta}\to{\Pi} \lambda_\alpha/_{J^{bd}_{\delta}}\right) =
\underset{\alpha < \delta}\to{(\sup}\, \lambda_\alpha)^+$.
\medskip
In particular, using Theorem 2.15 and the Strong Hypothesis, we see
that the following statement is true in $V^{\Bbb P}$.
\medskip
$(\ast)$  If $\mu$ is a singular cardinal such that $\mu = \sup
(CF(S)\cap\mu)$, then either\linebreak $cf(\mu)\in CF(S)$ or $\mu^+\in CF(S)$.
\medskip
This suggests that we might try to replace condition (1.4)(c) of
Theorem 1.3 by the following condition.
\medskip
$(1.4)(c)^\prime$  If $\mu$ is a singular
cardinal such that $\mu = \sup(C\cap\mu)$, then either\linebreak $cf(\mu)\in C$
or $\mu^+\in C$.
\medskip
However, Theorem 2.19 shows that this cannot be done.  For example,
Theorem 2.19 implies that if 
$$C=\{\aleph_1\} \cup \{\aleph_{\delta + 1}\big\vert \delta <
\omega_2, cf(\delta)=\omega\} \cup\{\aleph_{\omega_2+1}\},
$$
then there does not exist a c.c.c. notion of forcing $\Bbb P$ such
that $V^\Bbb P\vDash CF(S)=C$.
\remark{Remark 2.17}  The Strong Hypothesis is usually taken to be the
following apparently weaker statement.
$$
(2.18)  \text{ For all singular cardinals } \mu, pp(\mu)=\mu^+.
$$
(For the definition of $pp(\mu)$, see [Sh-400a].)  However,  Shelah
[Sh-420, 6.3 (1)] has shown that (2.16) and (2.18) are equivalent.
\endremark
\proclaim{Theorem 2.19 (The Strong Hypothesis)}
Let $\kappa$ be a regular uncountable cardinal, and suppose that
$\langle \lambda_\alpha|\alpha < \kappa\rangle$ is a strictly
increasing sequence of cardinals such that $\lambda_\alpha \in CF(S)$
for all $\alpha<\kappa$.  Suppose further that
\roster
\item"{(a)}"  $\kappa < \lambda_o$;
\item"{(b)}"  $E=\{\delta < \kappa\big\vert \lim \delta,
\underset{\alpha <\delta}\to{(\sup} \lambda_\alpha)^+ \notin CF(S)\}$
is a stationary subset of $\kappa$.
\endroster
Then $\kappa \in CF(S)$.
\endproclaim
\demo{Proof}  For each $\alpha < \kappa$, express
$S=\underset{i<\lambda_\alpha}\to{\bigcup} G^\alpha_i$ as the union
of a chain of $\lambda_\alpha$ proper subgroups.  For each $\delta \in
E$, let $\mu_\delta = \underset{\alpha <\delta}\to{\sup}\,
\lambda_\alpha$.  By the Strong Hypothesis, $tcf\left(\underset{\alpha
< \delta}\to{\Pi}
\lambda_\alpha/_{J^{bd}_\delta}\right)=\mu^+_\delta$. Let $\langle
f^\delta_\xi |\xi < \mu^+_\delta\rangle$ be a sequence in
$\underset{\alpha <\delta}\to{\Pi}\lambda_\alpha$ which witnesses that
$tcf\left(\underset{\alpha
<\delta}\to{\Pi}\lambda_\alpha/_{J^{bd}_\delta}\right)=\mu^+_\delta$.
For each $\xi < \mu^+_\delta$, let $H^\delta_\xi$ be the set of all
$g\in S$ such that $\delta\smallsetminus\{\alpha < \delta|g \in
G^\alpha_{f^\delta_\xi(\alpha)}\}\in J^{bd}_\delta$.  Once again, it is
easily checked that $\langle H^\delta_\xi|\xi<\mu^+_\delta\rangle$ is a
chain of subgroups such that
$S=\underset{\xi<\mu^+_\delta}\to{\bigcup} H^\delta_\xi$.  Since
$\mu^+_\delta \notin CF(S)$, there exists $\pi(\delta) < \mu^+_\delta$
such that $H^\delta_{\pi(\delta)}=S$.
\medskip
Since $\kappa < \lambda_o$, there exists $f\in\underset{\alpha<\kappa}\to{\Pi}\lambda_\alpha$ such that $f(\alpha) >
\sup\{f^\delta_{\pi(\delta)}(\alpha) |\alpha < \delta \in E\}$ for all
$\alpha < \kappa$.  Let $g \in S$.  Then for each $\delta \in E, g \in
H^\delta_{\pi(\delta)}$; and so there exists $\gamma(g,\delta) <\delta$
such that $g\in G^\alpha_{f^\delta_{\pi(\delta)}(\alpha)}\subseteq
G^\alpha_{f(\alpha)}$ for all $\gamma(g,\delta)\leq \alpha < \delta$.
By Fodor's Theorem, there exists an ordinal $\gamma(g)<\kappa$ and a
stationary subset $D$ of $E$ such that $\gamma(g,\delta)=\gamma(g)$ for
all $\delta \in D$.  This means that $g\in \bigcap\{
G^\alpha_{f(\alpha)} \gamma(g) \leq \alpha < \kappa\}$.
\medskip
For each $\gamma < \kappa$, let $\Gamma_\gamma =
\bigcap\{G^\alpha_{f(\alpha)} \gamma \leq \alpha < \kappa\}$.  Then
$\langle \Gamma_\gamma |\gamma <\kappa\rangle$ is a chain of subgroups
such that $S=\underset{\gamma < \kappa}\to{\bigcup}\Gamma_\gamma$.
Finally note that $\Gamma_\gamma\subseteq G^\gamma_{f(\gamma)}$, and so
$\Gamma_\gamma$ is a proper subgroup of $S$ for all $\gamma < \kappa$.
Thus $\kappa \in CF(S)$.
\enddemo
$\quad\hfill\qed$
\enddemo
\subhead 3.  The main theorem\endsubhead
In this section, we shall prove Thoerem 1.3.  Our notation generally
follows that of Kunen [K].  We shall only be using finite support
iterations.  An iteration of length $\alpha$ will be written as
$\langle \Bbb P_\beta, \tilde\Bbb Q_\gamma |\beta \leq \alpha, \gamma
< \alpha \rangle$, where $\Bbb P_\beta$ is the result of the first
$\beta$ stages of the iteration, and for each $\beta < \alpha$ there
is some $\Bbb P_\beta$-name $\tilde\Bbb Q_\beta$ such that
$$
\Vdash_{\Bbb P_\beta}\tilde\Bbb Q_\beta \text{ is a partial ordering}
$$
and $\Bbb P_{\beta +1}$ is isomorphic to $\Bbb P_\beta \ast \tilde\Bbb
Q_\beta$.  If $p\in \Bbb P_\alpha$, then $\text{ supt}(p)$ denotes the
support of $p$.
\medskip
There is one important difference between our notation and that of
Kunen.  Unlike Kunen, we shall {\it not} use $V^{\Bbb P}$ to denote
the class of $\Bbb P$-names for a notion of forcing $\Bbb P$.  Instead
we are using $V^\Bbb P$ to denote the generic extension, when we do
not wish to specify a particular generic filter $G\subseteq \Bbb P$.
Normally it would be harmless to use $V^\Bbb P$ in both of the above
senses, but there is a point in this section where this notational
ambiguity could be genuinely confusing.  Suppose that $\Bbb Q$ is an
arbitrary suborder of $\Bbb P$.  Then the class of $\Bbb Q$-names is
always a subclass of the class of $\Bbb P$-names.  (Of course, a $\Bbb
Q$-name $\tau$ might have very different properties when regarded as a
$\Bbb P$-name.  For example, it is possible that $\underset{\Bbb Q}\to{\Vdash}\tau$ is a
function, whilst $\underset{\Bbb P}\to{\nVdash}\tau$ is a function.)
However, we will not always have that $V^{\Bbb Q} \subseteq V^\Bbb P$;
where this means that $V[G\cap\Bbb Q] \subseteq V[G]$ for some
unspecified generic filter $G\subseteq\Bbb P$.
\definition{Definition 3.1} Let $\Bbb Q$ be a suborder of $\Bbb P$.
$\Bbb Q$ is a {\it complete} suborder of $\Bbb P$, written $\Bbb Q
\lessdot \Bbb P$, if the following two conditions hold.
\medskip
1.  If $q_1, q_2\in \Bbb Q$ and there exists $p\in \Bbb P$ such that
$p\leq q_1,q_2$, then there exists $r\in \Bbb Q$ such that $r\leq
q_1,q_2$.
\medskip
2.  For all $p\in \Bbb P$, there exists $q\in \Bbb Q$ such that
whenever $q^\prime \in \Bbb Q$ satisfies $q^\prime \leq q$, then
$q^\prime$ and $p$ are compatible in $\Bbb P$.  (We say that $q$ is a
{\it reduction} of $p$ to $\Bbb Q$.)
\medskip
It is wellknown that if $\Bbb Q \lessdot \Bbb P$, then $V^\Bbb
Q\subseteq V^\Bbb P$; and we shall only write $V^\Bbb Q\subseteq
V^\Bbb P$ when $\Bbb Q\lessdot \Bbb P$.
\medskip
We are now ready to explain the idea behind the proof of Theorem 1.3.
Let $V\vDash GCH$, and let $C$ be a set of regular uncountable
cardinals which contains a maximum element $\kappa$.  We seek a c.c.c.
$\Bbb P$ such that $V^\Bbb P \vDash 2^\omega = \kappa \wedge CF(S)=
C$.  The easiest part of our task is to ensure that $V^\Bbb P \vDash
C\subseteq CF(S)$.  We shall accomplish this by constructing $\Bbb P$
so that the following property holds for each $\lambda \in C$.
\enddefinition
$(3.2)_\lambda$  There exists a sequence $\langle\Bbb P^\lambda_\xi
|\xi < \lambda \rangle \in V$ of suborders of $\Bbb P$ such that 
\roster
\item"{(a)}"  if $\xi < \eta < \lambda$, then $\Bbb P^\lambda_\xi\lessdot\Bbb P^\lambda_\eta
\lessdot \Bbb P$;
\item"{(b)}"  for each $\pi \in$ Sym$(\omega)^{V^\Bbb P}$,
there exists $\xi < \lambda$ such that $\pi\in$ Sym$(\omega)^{V^{\Bbb
P^\lambda_\xi}}$;
\item"{(c)}"  for each $\xi< \lambda$, there exists $\pi\in$
Sym$(\omega)^{V^\Bbb P}\smallsetminus$ Sym$(\omega)^{V^{\Bbb
P^\lambda_\xi}}$.  
\endroster
\medskip
The harder part is to ensure that $V^{\Bbb P} \vDash CF(S) \subseteq
C$.  This includes the requirement that $(3.2)_\lambda$ fails for
every $\lambda\notin C$.  So, roughly speaking, we are seeking a
c.c.c. $\Bbb P$ which can be regarded as a ``kind of iteration'' of
length $\lambda$ precisely when $\lambda\in C$.  We shall use the
technique of Section 3 [Sh-288] to construct such a
notion of forcing $\Bbb P$.   
\definition{Definition 3.3}  Let $\langle a_i i<\alpha\rangle$ be a
sequence of subsets of $\alpha$.  We say that $b\subseteq \alpha$ is
{\it closed } for $\langle a_i|i<\alpha\rangle$ if $a_i \subseteq b$
for all $i\in b$.
\enddefinition
\definition{Definition 3.4}  Let $\Cal C$ be the class of all
sequences
$$\bar Q=\langle \Bbb P_i, \tilde \Bbb Q_j, a_ji\leq \alpha,
j<\alpha\rangle
$$
for some $\alpha$ which satisfy the following conditions.  (We say
that $\bar Q$ has length $\alpha$ and write $\alpha =$ lg $(\bar Q)$.)
\roster
\item"{(a)}"  $a_i\subseteq i$.
\item"{(b)}"  $a_i$ is closed for $\langle a_j |j< i\rangle$.
\item"{(c)}"  $\Bbb P_i$ is a notion of forcing and $\tilde \Bbb Q_j$
is a $\Bbb P_j$-name such that $\underset{\Bbb P_j}\to{\Vdash} \tilde
\Bbb Q_j$ is a c.c.c. partial order.
\item"{(d)}"  $\langle\Bbb P_i, \tilde \Bbb Q_j|i\leq \alpha , j
<\alpha\rangle$ is a finite support iteration.
\item"{(e)}"  For each $j<\alpha$, define the suborder $\Bbb
P^\ast_{a_j}$ of $\Bbb P_j$ inductively by 
\endroster
$$\Bbb P^\ast_{a_j} = \{p\in
\Bbb P_j|\text{ supt}(p) \subseteq a_j\text{ and } p(k)\text{ is a } \Bbb
P^\ast_{a_k}-\text{name for all } k\in \text{ supt}(p)\}.
$$

Then $\tilde \Bbb Q_j$ is a $\Bbb P^\ast_{a_j}$-name.  (At this stage, we
do not know whether $\Bbb P^\ast_{a_j}$ is a complete suborder of
$\Bbb P_j$.  It is for this reason that we are being careful with our
notation.  However, we shall soon see that $\Bbb P^\ast_{a_j} \lessdot
\Bbb P_j$, and then we can relax again.)
\enddefinition
\definition{Definition 3.5}  Let $\bar Q\in \Cal C$ be as above, so
that $\alpha =$ lg $(\bar Q)$.
\roster
\item"{(a)}"  We say that $b\subseteq \alpha$ is closed for $\bar Q$
if $b$ is closed for $\langle a_j|j<\alpha\rangle$.
\item"{(b)}"  If $b\subseteq \alpha$ is closed for $\bar Q$, then we
define $\Bbb P^\ast_b=\{p\in \Bbb P_\alpha \text{ supt}(p) \subseteq b$
and $p(k)$ is a $\Bbb P^\ast_{a_k}$-name for all $k\in \text{ supt}(p)\}$.
\endroster
\enddefinition
If $\beta < \alpha$, then we identify $\Bbb P_\beta$ with the
corresponding complete suborder of $\Bbb P_\alpha$ in the usual way.
If $b\subseteq \alpha$, then $p\upharpoonright b$ denotes the
$\alpha$-sequence defined by
$$\aligned 
(p\upharpoonright b)(\xi) &= p(\xi) \,\text{ if } \xi \in b\\
{ } &= \boldkey 1_{\tilde \Bbb Q_\xi}{} \,\,\text{ otherwise}
\endaligned
$$
\proclaim{Lemma 3.6}  Let $\bar Q\in \Cal C$ and let $\alpha = $lg
$(\bar Q)$.  Suppose that $b\subseteq c \subseteq \beta \leq \alpha$,
and that $b$ and $c$ are closed for $\bar Q$.
\medskip
(1)  $\beta$ is closed for $\bar Q$, and $\Bbb P_\beta = \Bbb
P^\ast_\beta$.
\medskip
(2)  If $p\in \Bbb P_\beta$ and $i\in \text{ supt}(p)$, then
$p\upharpoonright a_i \Vdash p(i) \in \tilde \Bbb Q_i$.
\medskip
(3)  Suppose that $p,q\in \Bbb P_\beta$ and $p\leq q$.  If $i\in
\text{ supt}(q)$, then $p\upharpoonright a_i \Vdash p(i) \leq q(i)$.
\medskip
(4)  If $p\in \Bbb P^\ast_c$, then $p\upharpoonright b \in \Bbb
P^\ast_b$.
\medskip
(5)  Suppose that $p\in \Bbb P^\ast_c, q \in \Bbb P^\ast_b$ and
$p\leq q$.  Then $p\upharpoonright b\leq q$.
\medskip
(6)  Suppose that $p\in \Bbb P^\ast_c, q\in \Bbb P_\beta$ and $p\leq q
\upharpoonright c$.
\medskip
Define the $\alpha$-sequence $r$ by
$$
\aligned
r(\xi) &= p(\xi) \text{ if } \xi \in c\\
{} &=q(\xi) \text  { otherwise}.
\endaligned
$$
Then $r\in \Bbb P_\beta$ and $r\leq p,q$.
\medskip
(7)  $\Bbb P^\ast_c\lessdot \Bbb P_\beta$.
\endproclaim
\demo{Proof} 
This is left as a straightforward exercise for the reader.
\enddemo
$\quad\hfill\qed$
\proclaim{Lemma 3.7}
Let $\bar Q \in \Cal C$ and let $\alpha =$ lg $(\bar Q)$.  Suppose
that $b\subset \alpha$ is closed for $\bar Q$ and that $i\in
\alpha\smallsetminus b$.
\medskip
(1)  $c=b\cup i$ and $c\cup\{i\}$ are closed for $\bar Q$.
\medskip
(2)  $\Bbb P^\ast_b\lessdot \Bbb P^\ast_c \lessdot\Bbb
P^\ast_{c\cup\{i\}}\lessdot \Bbb P_\alpha$.
\medskip
(3)  $\Bbb P^\ast_{c\cup\{i\}}$ is isomorphic to $\Bbb P^\ast_c \ast
\tilde \Bbb Q_i$.
\endproclaim
\demo{Proof}
Once again left to the reader.
\enddemo
$\quad\hfill\qed$
\medskip
Now we are ready to begin the proof of Theorem 1.3.  Suppose that
$V\vDash GCH$, and let $C$ be a set of regular uncountable cardinals
which satisfies the following conditions.
\medskip
(1.4)
\roster
\item"{(a)}"  $C$ contains a maximum element, say $\kappa$.
\item"{(b)}"  If $\mu$ is an inaccessible cardinal such that $\mu=\sup
(C\cap\mu)$, then $\mu \in C$.
\item"{(c)}"  If $\mu$ is a singular cardinal such that
$\mu=\sup(C\cap \mu)$, then $\mu^+\in C$.
\endroster

\definition{Definition 3.8}
\roster
\item"{(a)}"  $\Pi C$ denotes the set of all functions $f$ such that
dom $f = C$ and $f(\lambda)\in \lambda$ for all $\lambda \in C$.
\item"{(b)}"  $\Cal F_C$ is the set of all functions $f\in \Pi C$
which satisfy the following condition.
\item"{$(\ast)$}"  If $\mu$ is an inaccessible cardinal such that $\mu
=\sup (C\cap\mu)$, then there exists $\lambda < \mu$ such that $f(\theta)= 0$
 for all $\lambda \leq \theta \in C\cap \mu$.
\endroster
\enddefinition

\definition{Definition 3.9}
In $V$, we define a sequence 
$$
\langle \Bbb P_i, \tilde \Bbb Q_j, f_j i\leq \kappa, j <
\kappa\rangle
$$
such that the following conditions are satisfied.
\roster
\item"{(a)}"  $f_i \in \Cal F_C$.
\item"{(b)}"  Let $a_i=\{j< if_j \leq f_i\}$.  Then $\bar Q=\langle
\Bbb P_i,\tilde \Bbb Q_j, a_j|i\leq \kappa, j<\kappa\rangle \in \Cal
C$.
\item"{(c)}"  For each $f\in \Cal F_C$, there exists a cofinal set of
ordinals $j<\kappa$ such that $f_j=f$.
\item"{(d)}"  Suppose that $i<\kappa$ and that $\tilde \Bbb Q$ is a $\Bbb
P^\ast_{a_i}$-name with $|\tilde\Bbb Q|<\kappa$.  Then there exists
$i<j<\kappa$ such that 
\endroster
\medskip
(1)  $f_j=f_i$, and so $a_i\subseteq a_j$;
\medskip
(2)  if $\underset{\Bbb P_j}\to{\Vdash} \tilde \Bbb Q$ is c.c.c., then
$\tilde\Bbb Q_j=\tilde \Bbb Q$.
\enddefinition   
\medskip
We shall prove that $V^{\Bbb P_\kappa} \vDash CF(S) =
C$.  From now on, we shall work inside $V^{\Bbb P_\kappa}$.

\definition{Definition 3.10}  If $b\subseteq \kappa$ is closed for
$\bar Q$, then $S^b=$ Sym$(\omega)^{V^{\Bbb P^\ast_b}}$.
\medskip
First we shall show that $C\subseteq CF(S)$.  Fix some $\mu\in C$.  For
each $\xi < \mu$, let $b_\xi = \{i<\kappa|f_i(\mu) \leq \xi\}$.
Clearly $b_\xi$ is closed for $\bar Q$; and if $\xi<\eta <\kappa$,
then $b_\xi \subseteq b_\eta$.  Thus $\langle S^{b_\xi} \xi <
\mu\rangle$ is a chain of subgroups of $S$.
\enddefinition
\proclaim{Lemma 3.11}
For each $\xi < \mu, S^{b_\xi}$ is a proper subgroup of $S$. 
\endproclaim
\demo{Proof}
Let $\xi<\mu$ and let $i<\kappa$ satisfy $f_i(\mu)>\xi$.
Let $\Bbb Q$ be the partial order of finite injective functions
$q:\omega \rightarrow \omega$, and let $\tilde \Bbb Q$ be the
canonical $\Bbb P^\ast_{a_i}$-name for $\Bbb Q$.  Then there exists
$i<j<\kappa$ such that $f_j=f_i$ and $\tilde \Bbb Q_j=\tilde \Bbb Q$.
Clearly $j\notin b_\xi$.  Let $c=b_\xi \cup j$.  By Lemma 3.7, $\tilde\Bbb Q_j$ adjoins a permutation $\pi$ of $\omega$ such that
$\pi\notin V^{\Bbb P^\ast_c}$.  It follows that $\pi \notin
S^{b_\xi}$.
$\quad\hfill\qed$
\enddemo

\proclaim{Lemma 3.12}
$$
S=\underset{\xi<\mu}\to{\bigcup} S^{b_\xi}.
$$
\endproclaim
\demo{Proof}
Let $\pi\in S$.  Let $\tilde g$ be a nice $\Bbb
P^\ast_\kappa$-name for $\pi$. (Remember that $\Bbb P_\kappa=\Bbb
P^\ast_\kappa$.)  Thus there exist antichains $A_{\ell,m}$ of $\Bbb
P^\ast_\kappa$ for each $\langle\ell,m\rangle\in \omega\times \omega$
such that\newline $\tilde
g=\underset{\ell,m}\to{\bigcup}\{\langle\ell,m\rangle\}\times A_{\ell,m}$.
Let $\bigcup\{\text{ supt}(p) |p\in \underset{\ell,m}\to{\bigcup}
A_{\ell,m}\}=\{\alpha_k|k<\omega\}$.  Let\newline $\xi
=\sup\{ f_{\alpha_k}(\mu)|k<\omega\}$.  Then $p\in \Bbb P^\ast_{b_\xi}$
for each $p\in \underset{\ell,m}\to{\bigcup} A_{\ell,m}$, and so
$\tilde g$ is a nice $\Bbb P^\ast_{b_\xi}$-name.  Hence $\pi\in
S^{b_\xi}$.
\enddemo
$\quad\hfill\qed$
\medskip
This completes the proof of the following result.
\proclaim{Lemma 3.13}
If $\mu\in C$, then $\mu \in CF(S)$.\endproclaim
$\quad\hfill\qed$
\medskip
To complete the proof of Theorem 1.3, we must show that if $\mu\notin
C$, then\linebreak $\mu \notin CF(S)$.  We shall make use of the following easy observation.
\proclaim{Lemma 3.14}
Let $M\vDash ZFC$, and let $\langle g_\beta|\beta<\alpha\rangle\subseteq M$
be a generic sequence of elements of Sym$(\omega)$.  Let $\Bbb Q$ be
the partial order of finite injective functions
$q:\omega\rightarrow\omega$, and let $\pi\in M^\Bbb Q$ be the $\Bbb
Q$-generic permutation.  Then for all\newline $h\in$ Sym$(\omega)^M$, $\langle
g_\beta|\beta<\alpha\rangle\sphat\,  h\pi$ is generic.
\endproclaim
\demo{Proof}
For each finite subsequence $\beta_1<\cdots<\beta_n<\alpha$, the set
\newline $C(\alpha_1,\cdots,\alpha_n)=\{\phi
\in$ Sym$(\omega)|\langle g_{\alpha_1},\cdots,g_{\alpha_n}\rangle\sphat\, \phi$
is generic$\}$ is comeagre in Sym$(\omega)$.  Hence $h^{-1}C(\alpha_1,\cdots,\alpha_n)$ is also comeagre for each $h\in$
Sym$(\omega)$.  So for each\linebreak $h\in$
Sym$(\omega)^M,\pi \in h^{-1}C(\alpha_1,\cdots,\alpha_n)$.  The result
follows.
\enddemo
$\quad\hfill\qed$
\proclaim{Lemma 3.15}
Suppose that $\alpha <\kappa$ and that $\langle
g_\beta|\beta<\alpha\rangle$ is a generic sequence of elements of
Sym$(\omega)$.  If $H$ is any proper subgroups of Sym$(\omega)$, then
there exists a permutation $\phi\notin H$ such that $\langle
g_\beta|\beta <\alpha \rangle\sphat\, \phi$ is generic.
\endproclaim
\demo{Proof}
Let $h\in$ Sym$(\omega)\smallsetminus H$.  Then there exists
$i<\kappa$ such that $h,\langle g_\beta|\beta<\alpha\rangle \in
V^{\Bbb P_i}$.  There exists $i<j<\kappa$ such that $\tilde\Bbb Q_j$
is the canonical $\Bbb P^\ast_{a_j}$-name for the partial order $\Bbb
Q$ of finite injective functions $q:\omega\rightarrow\omega$.  By
Lemma 3.14, there exists a permutation $\pi\in V^{\Bbb P_{j+1}}$ such
that both $\langle g_\beta|\beta<\alpha\rangle\sphat\, \pi$ and $\langle
g_\beta|\beta<\alpha\rangle\sphat\,  h\pi$ are generic.  Clearly either
$\pi\notin H$ or $h\pi \notin H$.
\enddemo
$\quad\hfill\qed$
\medskip
Now fix some $\mu\notin C$, and suppose that $\mu\in CF(S)$.  It is
easily checked that $2^{\aleph_o}=\kappa$, and so we can suppose that
$\mu$ is a regular uncountable cardinal such that $\mu <\kappa$.
Express $S=\underset{\alpha <\mu}\to{\bigcup} G_\alpha$ as the union
of a chain of $\mu$ proper subgroups.  We can suppose that
Fin$(\omega)\leq G_o$.  Using Lemma 3.15, we can inductively construct
a generic sequence of elements of $S$
$$
\langle g^o_o, g^1_o,\cdots , g^o_\alpha, g^1_\alpha,\cdots
\rangle_{\alpha <\mu}
$$   
such that for each $\alpha <\mu$, there exists $\alpha\leq
\gamma_\alpha<\mu$ such that $g^o_\alpha \in G_{\gamma_\alpha}$ and
$g^1_\alpha\notin G_{\gamma_\alpha}$.
\proclaim{Lemma 3.16}
There exists a subset $X\in [\mu]^\mu$ and an ordinal $\xi <\kappa$
such that $\langle g^o_\alpha, g^1_\alpha \alpha \in X\rangle \in
V^{\Bbb P^\ast_{a_\xi}}$.
\endproclaim
\demo{Proof}
For each $\alpha <\mu$ and $\tau \in \{0,1\}$, let $\tilde
g^\tau_\alpha$ be a nice $\Bbb P^\ast_\kappa$-name for
$g^\tau_\alpha$.  Thus there exist antichains
$A^{\alpha,\tau}_{\ell,m}$ of $\Bbb P^\ast_\kappa$ for each $\langle
\ell,m\rangle \in \omega\times \omega$ such that 
$$
\tilde g^\tau_\alpha
= \underset{\ell,m}\to{\bigcup}\{\langle\ell,m\rangle\}\times
A^{\alpha,\tau}_{\ell,m}.
$$
For each $\alpha<\mu$, let
$\bigcup\{\text{ supt}(p) |p\in \underset{\ell,m}\to{\bigcup}A^{\alpha,o}_{\ell,m} \cup
\underset{\ell,m}\to{\bigcup}A^{\alpha,1}_{\ell,m}\} =\{\beta^\alpha_k
|k<\omega\}$.  Define $h_\alpha \in \Cal F_C$ by
$h_\alpha(\lambda)=\sup\{f_{\beta_k^\alpha}(\lambda) |k<\omega\}$ for
each $\lambda \in C$.
\medskip
It is easily checked that there are less than $\mu$ possibilities for
the restriction $h_\alpha\upharpoonright C\cap \mu$.  (This
calculation is the only point in the proof of Theorem 1.3 where we
make use of the hypothesis that $C$ satisfies conditions (1.4)(b) and
(1.4)(c).)  Hence there exists $X\in [\mu]^{\mu}$ such that
$h_\alpha\upharpoonright C \cap \mu=h_\beta\upharpoonright C\cap\mu$ for all $\alpha,\beta \in X$.
Define the function $f\in \Pi C$ by $f\upharpoonright C\cap \mu=
h_\alpha\upharpoonright C\cap \mu$, where $\alpha \in X$, and
$f(\lambda) =\sup\{h_\alpha(\lambda)|\alpha \in X\}$ for each $\lambda
\in C\smallsetminus\mu$.  Then it is easily checked that $f\in \Cal
F_C$; and clearly $f_{\beta^\alpha_k} \leq h_\alpha\leq f$ for all
$\alpha \in X$ and $k<\omega$.  Now choose
$\xi>\sup\{ \beta^\alpha_k|\alpha \in X,k< \omega\}$ such that $f_\xi
= f$.  If $\alpha\in X$ and $\tau\in\{0,1\}$, then $p\in \Bbb
P^\ast_{a_\xi}$ for each $p\in \underset{\ell,m}\to{\bigcup}
A^{\alpha,\tau}_{\ell,m}$; and hence $\tilde g^\tau_\alpha$ is a nice
$\Bbb P^\ast_{a_\xi}$-name.  It follows that $\langle g^o_\alpha,
g^1_\alpha \alpha\in X\rangle\in V^{\Bbb P^\ast_{a_\xi}}$.
\enddemo
$\quad\hfill\qed$
\medskip
By Lemma 2.13, there exists a $\sigma$-centred $\Bbb Q \in V^{\Bbb
P^\ast_{a_\xi}}$ such that 
$$
\underset{\Bbb Q}\to{\Vdash} \text{ There exists } \pi\in \text{
Sym}(\omega) \text{ such
that } \pi g^o_\alpha \pi^{-1}={}^\ast g^1_\alpha\text{  for all } \alpha \in X.
$$
Let $\tilde \Bbb Q$ be a $\Bbb P^\ast_{a_\xi}$-name for $\Bbb Q$.
Then there exists $\xi < \eta < \kappa$ such that $f_\eta = f_\xi$ and
$\tilde \Bbb Q_\eta = \tilde\Bbb Q$.  Hence there exists $\pi \in S$
such that $\pi g^o_\alpha\pi^{-1} ={}^\ast g^1_\alpha$ for all $\alpha
\in X$.  But this implies that $\pi\notin \underset{\alpha
<\mu}\to{\bigcup} G_\alpha$, which is a contradiction.  This completes
the proof of Theorem 1.3.
\medskip
By modifying the choice of the set $\Cal F_C$ of functions, we can
obtain some interesting variants of Theorem 1.3.  For example, the
following theorem shows that Theorem 2.2 cannot be proved in ZFC.
(Of course, it also shows that (1.4)(c) is not a necessary condition 
in Theorem 1.3.)
\proclaim{Theorem 3.17}  Suppose that $V\vDash GCH$ and that $\kappa >
\aleph_{\omega_1+1}$ is regular.  Let $C=\{\aleph_{\alpha +1} \alpha <
\omega_1\} \cup \{\kappa\}$.  Then there exists a c.c.c. notion of
forcing $\Bbb P$ such that $V^\Bbb P \vDash CF(S)=C$.
\endproclaim
\demo{Proof}  The proof is almost identical to that of Theorem 1.3.
The only change is that we use the set of functions
$$
\Cal F_C^\ast = \{f\in \Pi C|\text{ There exists }\alpha
<\omega_1\text{ such that } f(\aleph_{\beta +1})= 0\text{ for all }
\alpha \leq \beta < \omega_1\}
$$
in the definition of $\Bbb P_\kappa$.  This ensures that the counting
argument in the analogue of Lemma 3.16 goes through.
\enddemo
$\quad\hfill\qed$
\medskip
Using some more $pcf$ theory, we can prove the following result.
\proclaim{Theorem 3.18}
Suppose that $V$ satisfies the following statements.  
\roster
\item"{(a)}"  $2^{\aleph_n} = \aleph_{n+1}$ for all $n<\omega$.
\item"{(b)}"  $2^{\aleph_\omega} = \aleph_{\xi+1}$ for some $\omega
<\xi < \omega_1$.
\item"{(c)}"  $2^{\aleph_\eta} = \aleph_{\eta+1}$ for all $\eta \geq
\xi$.
\endroster
Let $T\in [\omega]^\omega$ and let $\kappa$ be a regular cardinal such
that $\kappa\geq \aleph_{\xi +1}$.  Let $C= pcf(\underset{n\in
T}\to{\Pi}\aleph_n)\cup\{\kappa\}$. Then there exists a c.c.c. notion
of forcing $\Bbb P$ such that $V^\Bbb P\vDash CF(S)=C$.
\endproclaim
\demo{Proof}  Again we argue as in the proof of Theorem 1.3.  This
time we use the set of functions, $\Cal F^\#_C = \underset{n\in
T}\to{\Pi}\aleph_n$, in the definition of $\Bbb P_\kappa$.  Examining
the proof of Lemma 3.16, we see that it is enough to prove that the
following statement holds for each regular uncountable $\mu\notin C$.
\medskip
$(3.19)_\mu$
$$
\aligned
&\text{ If }\langle h_\alpha|\alpha < \mu\rangle \text{ is a sequence
in } \underset{n\in T}\to{\Pi} \aleph_n, \text{ then there exists } X\in
[\mu]^\mu \\
&\text{ and an } f\in\underset{n\in T}\to{\Pi} \aleph_n
\text{ such that} h_\alpha \leq f \text{ for all } \alpha \in X.
\endaligned
$$
\medskip
This is easy if $\mu<\aleph_\omega$.  If $\mu> \aleph_\omega$, then
$(3.19)_\mu$ is a consequence of the following result.
\enddemo
\proclaim{Theorem 3.20}  Let $\{\lambda_i|i\in I\}$ be a set of
regular cardinals such that\newline $\min\{\lambda_i|i\in I\}>|I|$.  Let $\mu$
be a regular cardinal such that $\mu>2^{|I|}$ and\newline $\mu \notin pcf
(\underset{i\in I}\to{\Pi} \lambda_i)$.  If $\langle h_\alpha|\alpha
<\mu\rangle$ is a sequence in $\underset{i\in I}\to{\Pi}\lambda_i$,
then there exists $X\in [\mu]^\mu$ and $f\in \underset{i\in
I}\to{\Pi}\lambda_i$ such that $h_\alpha \leq f$ for all $\alpha \in
X$.
\endproclaim
\demo{Proof}
This is included in the proof of [Sh-g, II 3.1].
(More information on this topic is given in [Sh-513, Section 5].
Also [Sh-430, 6.6D] gives even more information under the hypothesis
that $2^{|I|} < \min\{\lambda_i|i\in I\}$.)  Alternatively, argue as in the proof
of [BM,7.11].
\enddemo
$\quad\hfill\qed$
\medskip
It is known that, assuming the consistency of a suitable large
cardinal hypothesis, for each $\omega<\xi<\omega_1$ there exists a
universe which satisfies the hypotheses of Theorem 3.18.  (See [GM].)
Thus the following result shows that Theorem 1.2 cannot be
substantially improved in $ZFC$.
\proclaim{Corollary 3.21}  Suppose that $V$ satisfies the hypotheses
of Theorem 3.18 with respect to some $\omega<\xi<\omega_1$.  Then for
each $\omega \leq \alpha \leq\xi $ and $\kappa\geq \aleph_{\xi +1}$, there
exists a set $T\in [\omega]^\omega$ and a c.c.c. notion of forcing
$\Bbb P$ such that 
$$
V^\Bbb P \vDash CF(S) = \{\aleph_n n\in T\}\cup\{\aleph_{\alpha
+1}\}\cup\{\kappa\}.
$$
In particular, if $\omega <\alpha \leq \xi$, then 
$$
V^\Bbb P\vDash\aleph_{\omega + 1}\notin CF(S).
$$
\endproclaim
\demo{Proof}
With the above hypotheses, [Sh-g,VIII] implies that there exists $T\in
[\omega]^\omega$ such that $tcf\left(\underset{n\in T}\to{\Pi}
\aleph_n/_{J^{bd}_\omega}\right) = \aleph_{\alpha +1}$.  It follows
that\newline $pcf (\underset{n\in T}\to{\Pi}\aleph_n)=\{\aleph_n|n\in
T\}\cup\{\aleph_{\alpha +1}\}$.  So the result is a consequence of
Theorem 3.18.
\enddemo
$\quad\hfill\qed$
\medskip
Finally we shall show that (1.4)(a) is not a necessary condition in Theorem 1.3,
and that $2^{\aleph_o}$ cannot be bounded in terms of the set $CF(S)$.
\proclaim{Theorem 3.22} Suppose that $V\vDash GCH$ and that 
$C=\{\aleph_{\alpha +1} \alpha < \omega_1\}$. If $\kappa$ is any singular
cardinal such that $cf(\kappa)\in C$, then there exists a $c.c.c$ notion of
forcing $\Bbb P$ such that $V^\Bbb P \vDash CF(S)=C \text{ and } 
2^{\aleph_o} = \kappa$.
\endproclaim
\demo{Proof}
Let $\kappa$ be a singular cardinal such that $cf(\kappa) \in C$. Let
$\langle \lambda_\beta|\beta < cf(\kappa)\rangle$ be a strictly increasing
sequence of regular cardinals such that $\lambda_0 = \aleph_{\omega_1+1}$
and $\underset{\beta < cf(\kappa)}\to{\sup}\, \lambda_\beta =\kappa$. Let
$$
\Cal F_C^\ast = \{f\in \Pi C|\text{ There exists }\alpha
<\omega_1\text{ such that } f(\aleph_{\beta +1})= 0\text{ for all }
\alpha \leq \beta < \omega_1\}.
$$ 
In V, we define a sequence $\langle \Bbb P_i,\tilde \Bbb Q_j, f_j i\leq \kappa,
j < \kappa\rangle$ such that the following conditions are satisfied.
\roster
\item"{(a)}" $f_i \in \Cal F_C^\ast$.
\item"{(b)}" Let $a_i=\{j< if_j \leq f_i\}$. Then $\bar Q=\langle
\Bbb P_i,\tilde \Bbb Q_j, a_j|i\leq \kappa, j<\kappa\rangle \in \Cal
C$.
\item"{(c)}" For each $f\in \Cal F_C^\ast$ and $\beta < cf(\kappa)$, there exists
a cofinal set of ordinals $j < \lambda_\beta$ such that $f_j=f$.
\item"{(d)}" Suppose that $\beta < cf(\kappa)$, $i < \lambda_\beta$ and that
$\tilde \Bbb Q$ is a $\Bbb P^\ast_{a_i}$-name with $|\tilde\Bbb Q|<\lambda_\beta$.
Then there exists $i<j<\lambda_\beta$ such that
\endroster
\medskip
(1) $f_j=f_i$, and so $a_i\subseteq a_j$;
\medskip
(2) if $\underset{\Bbb P_j}\to{\Vdash} \tilde \Bbb Q$ is c.c.c., then
$\tilde\Bbb Q_j=\tilde \Bbb Q$.
\medskip
Clearly $V^{\Bbb P_\kappa} \vDash 2^{\aleph_0} = \kappa$. Arguing as in the 
proof of Lemma 3.13, we see that $V^{\Bbb P_\kappa} \vDash C\subseteq CF(S)$.
>From now on, we shall work inside $V^{\Bbb P_\kappa}$.
Let $\mu$ be a regular cardinal such that
 $\aleph_{\omega_1+1} \leq \mu < \kappa$.
 Suppose that we can express
$S=\underset{\alpha <\mu}\to{\bigcup} G_\alpha$ as the union of a chain of
$\mu$ proper subgroups. For each $\alpha < \mu$, choose an element
$h_\alpha \in G\smallsetminus G_\alpha$. Then there exists a subset 
$I\in [\mu]^\mu$ and an ordinal $\beta < cf(\kappa)$ such that
$\langle h_\alpha \alpha \in I\rangle \in V^{\Bbb P_{\lambda_\beta}}$ and
$\mu\leq \lambda_\beta$. In $V^{\Bbb P_\kappa}$, we can inductively construct
a generic sequence of elements of $S$
$$
\langle g^0_0, g^1_0,\cdots , g^0_\alpha, g^1_\alpha,\cdots
\rangle_{\alpha < \mu}
$$
such that for each $\alpha <\mu$
\roster
\item"{(1)}" there exists $\alpha\leq \gamma_\alpha < \mu$ such that
$g^0_\alpha \in G_{\gamma_\alpha}$ and $g^1_\alpha\notin G_{\gamma_\alpha}$; and
\item"{(2)}" there exists $\lambda_\beta \leq i_\alpha < \lambda_{\beta +1}$ such
that $\langle g^0_\delta, g^1_\delta |\delta < \alpha \rangle
\subseteq V^{\Bbb P_{i_\alpha}}$.
\endroster
\medskip
For suppose that $\langle g^0_\delta, g^1_\delta |\delta < \alpha\rangle$ has
been defined. By Lemma 3.14, there exists $i_\alpha < j < \lambda_{\beta +1}$
and $g^0_\alpha \in V^{\Bbb P_j}$ such that 
$\langle g^0_\delta, g^1_\delta |\delta < \alpha\rangle\sphat\, g^0_\alpha$
is generic. Choose $\gamma_\alpha \in I$ such that $\alpha\leq \gamma_\alpha < \mu$
and $g^0_\alpha \in G_{\gamma_\alpha}$. By a second application of Lemma 3.14, there
exists $j < i_{\alpha +1} < \lambda_{\beta +1}$ and $\pi \in V^{\Bbb P_{i_{\alpha +1}}}$
such that both $\langle g^0_\delta, g^1_\delta |\delta < \alpha \rangle\sphat\, 
g^0_\alpha{ \sphat}\, \pi$ and
$\langle g^0_\delta, g^1_\delta |\delta < \alpha\rangle\sphat\, {g^0_\alpha} \sphat\,
h_{\gamma_\alpha}\pi$ are generic. Clearly either $\pi \notin G_{\gamma_\alpha}$
or $h_{\gamma_\alpha}\pi \notin G_{\gamma_\alpha}$. Hence we can also find a
suitable $g^1_\alpha$.
\medskip
There exists a subset $J \in [\mu]^\mu$ and an ordinal $\delta < cf(\kappa)$
such that\newline $\langle g^0_\alpha, g^1_\alpha |\alpha \in J\rangle \in V^{\Bbb P_{\lambda_\delta}}$
and $\mu\leq \lambda_\delta$. Arguing as in the proofs of Theorems 1.3 and 3.17, there
exists $\pi \in V^{\Bbb P_{\lambda_{\delta +1}}}$ such that
$\pi g^0_\alpha \pi^{-1}={}^\ast g^1_\alpha$ for all $\alpha \in J$. This is a contradiction.
\enddemo
$\quad\hfill\qed$

\bigskip

\medskip

\noindent{MATHEMATICS DEPARTMENT}
\smallskip
\noindent{BILKENT UNIVERSITY}
\smallskip
\noindent{ANKARA}
\smallskip
\noindent{TURKEY}

\bigskip

\medskip

\noindent MATHEMATICS DEPARTMENT
\smallskip
\noindent THE HEBREW UNIVERSITY
\smallskip
\noindent JERUSALEM
\smallskip
\noindent ISRAEL

\bigskip

\medskip

\noindent{MATHEMATICS DEPARTMENT}
\smallskip
\noindent{RUTGERS UNIVERSITY}
\smallskip
\noindent{NEW BRUNSWICK, NEW JERSEY}
\smallskip
\noindent{USA}

\end